# Supplementary Private Tutoring and Mathematical Achievements in Higher Education:

## An Empirical Study on Linear Algebra


Xuefei Lin[1,2], Guangyu Xu[1,2], Lei Peiyao[3], Bin Xiong[1,2]

1. School of Mathematical Sciences, East China Normal University,

2. Shanghai Key Laboratory of Pure Mathematics and Mathematical Practice

3. Chengdu Academy of Education Sciences, Chengdu, China



**Funding**

This research was supported by Shanghai Key Laboratory of Pure Mathematics and Mathematical Practice. (Fund No.: 18DZ2271000)



**Correspondence**

Xuefei Lin    Linxuefei@stu.ecnu.edu.cn

Guangyu Xu 52205500036@stu.ecnu.edu.cn

Peiyao Lei    leipyy@163.com

Bin Xiong    bxiong@math.ecnu.edu.cn





**Abstract**

The present article is an empirical study that investigates the learning situation of linear algebra. Research was performed among 60 science and engineering students from different universities in Zhejiang, Jiangsu, Hubei, and Shandong who were taught by three mathematics teachers in the course of a 3-day short-term training consisting of 24 classes; following this, the effectiveness of short-term training for mathematics curriculum learning in higher education was evaluated based on the IRT measurement model. The results indicate that regardless of the difficulty of the test, short-term training significantly enhances students' mastery of linear algebra knowledge points, and students can accurately perceive their mathematical achievements; the effect observed is not consistent with that of private tutoring (PT) in primary and secondary schools.

**Supplementary Private Tutoring and Mathematical Achievements in Higher Education:**

**An Empirical Study on Linear Algebra**

Linear algebra is a captivating field that is often introduced to students in the first two years of college. This field emphasizes the interplay between definitions and proofs (Axler, 1997). It diverges from students' previous academic exposure by introducing an array of novel concepts (Strang, 2022). Principal among these are linear independence, linear transformations, vector space, invertibility, and rank. These notions must be employed deliberately and connected to substantiate theorems or solve problems. The pervasive utilization of definitions coupled with theorem dependencies often lends the discipline an abstract and theoretical character (Lang, 2012). While such a mode of teaching appeals to experts, it may be off-putting for many students. Given the challenges some students face in assimilating the subject matter, considerable attention was extended to this topic by ICMI-13, which culminated in the publication of the book Challenges and Strategies in the Teaching of Linear Algebra, aiming to address these issues.

Certain representative concerns in linear algebra have piqued the interest of several researchers (Dias & Artigue, 1995; Hillel, 2000). Numerous inquiries have unveiled various facets of learning regarding the abstract constituents of linear algebra using diverse theoretical perspectives, including vector space (Dorier, 1995), linear transformation (Sierpinska et al., 1999), linear combination, and linear independence (Dogan-Dunlap, 2018). Pivoting to a specific mathematical notion, linear transformations serve as a cornerstone in linear algebra. These functions, which preserve linearity between vector spaces, permit an amplified study of vector space structures through the investigation of their properties.



However, this concept inherently poses significant challenges for students in their understanding and for educators in their teaching endeavors (Dorier & Sierpinska, 2001).

The aforementioned studies strive to elucidate the complexities associated with the comprehension of linear algebra's knowledge points. Viewed holistically, confronting the challenges of linear algebra is itself a daunting task. For instance, Ferryansyah conducted a study using written tests and interviews administered to 35 students. The findings revealed that the difficulty students experienced in understanding linear algebra was exceptionally high; 88.63% were not adept at symbol or notation representation, and 88.11% struggled with the implementation of mathematical symbols, notations, or logical reasoning. Moreover, 88.38% found it challenging to comprehend symbols or notations through logical reasoning. Finally, 91.77% had difficulty verifying the correct application of mathematical symbols, notations or concepts and employing logical reasoning (Ferryansyah et al., 2018).

The complexity of linear algebra often leads to students failing their examinations. In mainland China, proficiency in linear algebra is a prerequisite for postgraduate applications, which compels many students to seek supplementary learning tools to attain their objectives. Private supplementary tutoring (PT), colloquially known within the academic research community as shadow education, mirrors the formal school curriculum (Baker et al., 2001; Bray & Kwok, 2003).

According to the Quality Monitoring Report on China's Compulsory Education released by the National Center of Education Quality (Zhang et al., 2021), 23.4% of eighth graders in China participated in after-school mathematics classes. Supplementary private tutoring (PT) has become a prevalent practice in many regions, including Europe and the



USA, in addition to East Asia, and has been increasingly discussed by researchers worldwide (Bray, 2014; Davies, 2004; Mischo & Haag, 2002; Zhang et al., 2021). This widespread social phenomenon has significant impacts on the daily lives of those involved in education (Liu & Bray, 2018) and has opened up new research areas (Bray, 2014; Davies, 2004) that extend beyond the focus on regular school or family education (He et al., 2021). PT is a basic strategy that may be adopted by parents, students, teachers, principals, and school systems when students fail to meet academic standards (Bray, 2014; Wang & Guo, 2017; Zeng & Zhou, 2012).

In the Canadian educational landscape, certain coaching establishments, colloquially termed "cram centers," posit—sans empirical corroboration—that sustained training on advanced cognitive skills can yield significant academic dividends (Davies, 2004). Intriguingly, this assertion finds resonance in the context of contemporary China. However, many private tutoring (PT) institutions use these claims in their aggressive marketing strategies, despite a conspicuous absence of systematic evidence attesting to their programmatic effectiveness. A pivotal inquiry within the scholarship on this subject is the actual impact of PT in mathematics on students' longitudinal academic trajectories (He et al., 2021; Zhang et al., 2021). Given the substantial financial and time commitments families allocate to PT in mathematics, they are unequivocally justified in seeking tangible returns on these investments, particularly in the form of enhanced academic outcomes in mathematics (Bray, 2014).

Private tutoring (PT) is prevalent not only in primary and secondary schools but also in universities. Thus, this article aims to investigate the prevalence of PT in the context of



linear algebra in higher education.

**Research Questions**

RQ1: Does short-term training contribute to noticeable improvement in students' achievement in linear algebra at a higher education level? Further, does short-term training enhance students' problem-solving abilities? In which ways does it do so?

RQ2: Are students able to accurately identify their own variations within the context of linear algebra remediation?

## Literature Review

**Difficulties in Linear Algebra Materials**

*Determinants*

Donevska-Todorova (2014) classified determinant thought processes into three descriptive modes: comprehensive geometric thinking (or geometric language), analytical arithmetic thinking (or arithmetic language), and analytical structure thinking (or algebraic language). In the comprehensive geometric mode, the determinant is interpreted as the oriented volume and area of a parallelepiped and a parallelogram, individually expanded by vectors, while in the analytical arithmetic mode, the determinant is considered to be the sum of permutations. Finally, in the analytical structure mode, the determinant is defined as a function that fulfills three axioms: multilinear form, a norm, and zero value for the determinant of two identical rows in the matrix.

Additionally, Aygor and Ozdag (2012) researched the misconceptions held by 60 undergraduates in solving determinant and matrix problems. Their findings showed numerous misunderstandings, primarily involving confusion between a matrix and its determinant. For



example, many students make the following conclusion: if $A = kB$, then $\det A = k \det B$. This confusion can often lead to poor academic performance in mathematics. Remarkably, similar errors have also been found among teachers specializing in linear algebra computation. A study in Zimbabwe involving 116 teachers identified 34 participants who made similar mistakes (Kazunga & Bansilal, 2018).

### *Matrix Multiplication*

Cook et al. (2018) conducted a comprehensive analysis of various instructional materials related to matrix multiplication. Renowned mathematicians highlight four pivotal roles of matrix multiplication, specifically, its functions as a system of linear equations, a linear space, a linear transformation, and an elementary row transformation. Cook underscored that the introduction of matrix multiplication invariably influences students' comprehension. Further enriching the theoretical framework proposed by Harel (1987), Cook incorporated aspects of computational efficiency and seamlessly integrated concepts of isomorphism, delay, abstraction, analogy, and computational efficacy. This refined model has proven instrumental in advancing the understanding of matrix multiplication.

### *Systems and Solutions*

Panizza et al. (1999) asked students to independently construct a 2×2 linear system and had them solve it. While students could adeptly determine the system's solution, they struggled with discerning whether that solution also satisfied specific individual equations within the system. The two-variable equation was not perceived as a defining parameter for an infinite set of logarithms. From these observations, Panizza inferred that students grapple with articulating the essence of equations, potentially impacting their understanding of



systems formed by these equations.

Nonetheless, in linear algebra instruction, the emphasis should not be exclusively on the concept of solutions. Instead, a significant focus should be directed toward methods for solving equations, such as the Gauss algorithm. The previously discussed APOS approach presents challenges in assimilating such algorithms. This is primarily because students find it daunting to envision outcomes without going through the concrete steps of problem solving (DeVries & Arnon, 2004).

### *Vectors*

Linear independence has emerged as the focal point in several notable studies, notably in the research by Stewart et al. (2019). An intriguing study by Ertekin et al. (2010) explored 144 primary and high school mathematics teacher trainees' grasp of the concepts of linear independence and dependence. In contrast with the findings of Beltrán-Meneu et al. (2016), Ertekin's team discovered that students displayed greater adeptness in algebraic facets as opposed to geometric ones. Dogan-Dunlap (2010, p. 2141), venturing into the nexus between student comprehension and graphical representations pertaining to linear independence, posited that "the juxtaposition of geometric representations alongside algebraic and arithmetic modes seemingly empowers learners to fluidly navigate the multifaceted representational dimensions of a concept". In parallel, other scholars have devised instructional sequences to bolster students' insights about span and linear independence, meticulously intertwining algebraic and geometric representations (Wawro et al., 2013).

Dogan-Dunlap (2010) also framed her study using Sierpinska's (2000) modes of thinking. Students in the Dogan-Dunlap study completed a homework assignment using an



interactive geometric web module in which they solved problems involving linear independence and dependence of vectors in R3. Student homework solutions revealed 17 categories of student thinking modes, 11 of which were geometric. However, students tended to use the 6 algebraic or arithmetic thinking modes in explaining whether or not sets were linearly independent and in developing conjectures about linearly independent sets of vectors in R3. Dogan-Dunlap (2018) continued this line of research with interactive web modules and refined her previous categories into 9 student thinking modes, each of which were further classified into categories of those allowing for 1, 2 or all 3 of Sierpinska's modes of thinking. In this study, three different linear algebra classes had different levels of access to the web modules, as follows: no access, homework only, or homework and in-class demonstration. Dogan found that students were more likely to reason geometrically in follow-up interviews about linear independence if they had used the interactive geometric modules outside class and even more so if those modules had also been used by their teacher in class.

The observed student misconceptions in linear algebra, such as the inability of students to correctly recognize the determinant $\det kA \neq k \det A$ and other related conclusions, as well as the difficulty of component algebra models, were drawn upon to create the material for the test questions.

### *The Roles and Functions of PT in Mathematics Learning*

The objective of private tutoring (PT) is to amplify the outcomes of traditional educational systems and broaden their curricular scope (Zeng & Zhou, 2012). In this form of supplementary education, students may either advance their learning or reinforce previously acquired school knowledge. Typically, PT primarily targets students at both ends of the



mathematical proficiency spectrum.

For the academically adept, tutoring serves as a platform to further hone their strengths or, at a minimum, sustain their academic standing. In contrast, for those who are struggling, PT provides a valuable avenue to bridge the knowledge gap and align with their peers (Bray & Kobakhidze, 2014). An additional objective of PT often pertains to specialized academic competitions. A case in point is the emphasis on elevating performance levels in mathematical Olympiads in China (Wang & Guo, 2017). Moreover, there is a growing emphasis on courses tailored for autonomous admissions in local universities and preparatory modules designed for international academic pursuits, such as advanced placement programs, which have gained significant traction in recent years.

### *Effects of PT on Students' Mathematics Learning Performance*

Central to the global PT (private tutoring) research paradigm, the field has been dominated by studies investigating its effectiveness. These seminal investigations not only elucidate the superior performance of Asian students in international assessment frameworks but also serve as critical guides for parental educational choices. Building upon an expansive theoretical framework, future research should focus on performing a thorough evaluation of PT efficacy. Such assessments demand rigorously constructed data analytics processes paired with representative samples to yield reliable outcomes (Ren et al., 2022). To ascertain the extent to which PT augments students' cognitive faculties and enhances their inherent problem-solving proficiencies, item response theory (IRT) emerges as a potent analytical tool.

The overall effectiveness of PT should be tested, and the instructional quality should be analyzed to explain the corresponding results. Furthermore, to counteract the hyperbolic



and, at times, misleading promotional claims surrounding PT, robust empirical studies are indispensable. Such scholarly endeavors will not only guide policy determinations by educational authorities but also help families refine their educational strategies.

***Item Response Theory***

Item response theory (IRT) stands as a foundational pillar in the domain of psychometrics, offering a sophisticated framework that delineates the nuanced interplay between individuals' latent traits—which are often their abilities or proficiencies—and their consequent responses to assessment items (Hambleton et al., 1991).

Historically, classical test theory (CTT) was the dominant form of assessment interpretation. However, its principal focus was on aggregate test-level indices, often proving inadequate for detailed item-level scrutiny. IRT emerges as a transformative alternative, as it emphasizes the attributes of individual test items and facilitates inferences regarding test-takers and items on a singular, unified scale (Lord, 2012).

One of IRT's salient strengths lies in its invariance property: it allows the constancy of item parameters across disparate populations, contingent upon model fit and unbiased trait measurement (Embretson & Reise, 2013). Such invariance enables a range of applications, including equating scores across different test forms and facilitating the advancement of computer-adaptive testing (CAT) methodologies (Wainer et al., 2000).

IRT's superiority in test development and validation spheres has been nothing short of revolutionary. Its potent capability to augment measurement precision, together with advanced tools for differential item functioning (DIF) analyses, ensures that items retain consistent behavior across diverse demographic subgroups, thereby addressing potential



biases (Zumbo, 2007).

Modern educational landscapes are increasingly globalized, with an escalating demand for standardized, comparable assessments. Here, IRT's relevance and applicability stand out prominently. Notably, contemporary advancements in computational techniques and specialized software tools have democratized IRT's adoption, underscoring its indispensable role in contemporary psychometrics (DeMars, 2010).

The sheer depth and granularity of insights provided by IRT make it an invaluable asset for researchers and practitioners alike. As assessments and evaluations continue to evolve, especially in an era of digital proliferation, the dynamic interplay of theory, methodology, and practice facilitated by IRT will undoubtedly remain at the forefront of educational measurement (Samejima, 1997).

### Theoretical Framework

For the first research question, on the basis of the textbook Linear Algebra (2021) and previous articles by Aygor and Ozdag (2012), Kazunga & Bansilal (2018), Cook (2018), Panizza (1999), DeVries & Arnon (2004), Ertekin (2010) and others mentioned articles, the linear algebra studied by undergraduate students has been categorized as follows according to the knowledge points, and the results are shown in Table 1.

**Table 1**

*Linear Algebra Breakdown of Knowledge Points*

| Topic | Knowledge points | Pre-test | Post-test |
|---|---|---|---|
| Determinant | The concept and basic properties of determinants | 1 | 21 |
| | Adjugate matrix | 2 | 17 |
| | Multiplication of matrices | 5 | 22 |
| | The power of a square matrix | 5 | 22 |
| Matrix | The concept and properties of inverse matrices | 6 | 18 |
| | Rank of matrix | 7 | 19 |
| | Elementary matrix | 11 | 27 |
| | | 3 | 23 |
| Vector | Linear combination and linear representation of vectors | 14 | 30 |
| | Linear correlation and linear independence of vector groups | 3 | 23 |
| | The maximum linearly independent group of vector groups | 14 | 30 |
| | Rank of vectors | 14 | 30 |
| System of linear equations | The properties of solutions to linear equations and the structure of solutions | 8 | 20 |
| | | | 29 |
| | The basic solution system and general solution of homogeneous linear equations | 13 | 29 |
| | General solutions for non homogeneous linear systems of equations | 13 | 29 |
| Eigenvalues and eigenvectors | The concept and properties of similar matrix | 9 | 25 |
| | The concept and properties of eigenvalues and eigenvectors of matrices | 10 | 26 |
| | Sufficient and necessary conditions for matrix to be similar diagonalized and similar diagonal matrix | 15 | 31 |
| Quadratic form | Positive definiteness of quadratic form and its matrix | 4 | 24 |
| | Inertia theorem | 12 | 28 |
| | Standard form and normal form of quadratic form | 16 | 32 |
| | Orthogonal transformation and matching method | 16 | 32 |

Formal and out-of-school mathematics instruction should be included when discussing potential factors that promote student learning in mathematics or analysing differences in student achievement in mathematics (Wang, and Guo, 2017; Wang et al ., 2019; Zhang et al, 2021). However, their effectiveness on student learning may vary depending on several characteristics of PT (Bray, 2014), such as scheduling (during school terms or holidays), which has not been examined by previous researchers.PT adds external time to student learning, and the amount of time that students need to learn depends not only on the characteristics of the students (Hiebert, and Grouws 2007), but also on the quality of teaching (Guill et al ., 2020; Zhang et al, 2021).Zhang (2021) used Carroll's (1970) model of school learning and Dunkin and Biddle's (1974) model of teacher effects to organise the aspects of physical education teaching quality that affect students' performance in mathematics as shown in Table 2. Several predictor or process variables discussed in previous studies as well



as some newly developed variables to explain the effect of PT on the product variable (school mathematics achievement) were added to this framework. Based on the above discussion, the theoretical framework in Table 2 below was obtained based on Carroll's (1970) model of school learning.

**Table 2**
*Theoretical Framework of RQ2*

| Dimension | Sub-dimension |
| --- | --- |
| Product variables | Post-mathematics achievement |
| Process variables | Instructional content |
| Presage variables | Teacher qualifications; Time point; Forms of instructional organization |
| Context variables | Premathematics achievement; Self-expectation of education; Family SES; Demographic characteristics |
| Time spent | Hours spent on PT |
| Intensity | |

## Methodology

### Participants and Procedures

#### *Prior to the Test*

Since the students' trait levels may change before and after the teaching intervention, it is not appropriate to combine the pre-test and post-test scores for IRT analysis. Therefore, a preliminary test is needed before the formal test to estimate parameters such as the difficulty, discrimination, and guessing of the test items. In order to ensure the quality of the test papers, more than 700 volunteers were recruited to take the test from 2021 to 2022, respectively, at the schools of the volunteers where the formal test was implemented. For this test to run



smoothly, the researcher collaborated with a number of local organizations that provide training for postgraduate entrance exams in the form of admissions counselling and disrupted and fused the two sets of test questions into a single set of questions. In the end, a total of 659 questionnaires were collected from the pilot test to calculate the item parameters, including difficulty, discrimination, and guessing, using IRT.

### The Fit of Model

First, the fit of the 1PL, 2PL, and 3PL models was evaluated. Based on AIC, SABIC, BIC, and HQ, 2PL is the optimal model because it has the lowest penalized information criterion values are shown in table 3. Although the log-likelihood of 3PL is higher, indicating a slightly better fit, due to its higher model complexity (with higher AIC and BIC) and the chi-square test showing no significant improvement, mod1 remains the most appropriate model overall.

**Table 3**
*The Fit of Model*

|       | AIC      | SABIC    | HQ       | BIC      | logLik   |
|-------|----------|----------|----------|----------|----------|
| 1PL   | 24647.93 | 24732.14 | 24759.34 | 24935.34 | 12259.97 |
| 2PL   | 24620.98 | 24664.40 | 24678.42 | 24769.17 | 12277.49 |
| 3PL   | 24685.39 | 24811.69 | 24852.50 | 25116.50 | 12246.69 |

The p-value for the 2PL model is 0.374, the RMSEA is very low (0.005), and the SRMSR is also small (0.032). Additionally, both TLI and CFI are close to 1, indicating that the model fits well and explains the data effectively. Therefore, the model's fit is highly satisfactory, with little need for further improvement.

### Unidimensionality Test

The unidimensionality assumption is a key assumption in Item Response Theory



(IRT). It assumes that all test items in a given test measure the same latent trait or ability, and that this latent trait is singular. Exploratory factor analysis was conducted and a scree plot was generated, as shown in the figure1 below. From the data in the figure, it can be observed that the ratio of the first eigenvalue to the second eigenvalue is greater than 3, indicating that the unidimensionality assumption is satisfied.

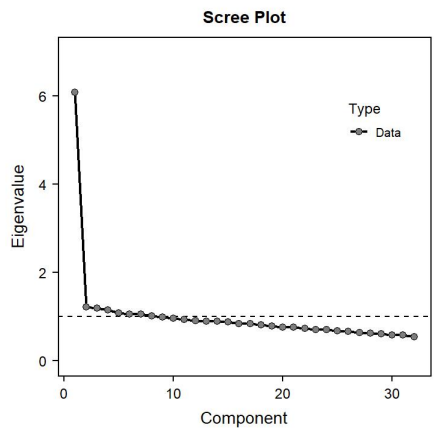

**Figure 1**
*Scree Plot of Factor Analysis*

Exploratory factor analysis was conducted and a scree plot was generated, as shown in the figure below. From the data in the figure, it can be observed that the ratio of the first eigenvalue to the second eigenvalue is greater than 3, indicating that the unidimensionality assumption is satisfied.

### Local independence test

Local independence is a common concept in statistics, particularly in factor analysis and item response theory. The core idea is that, after controlling for certain factors, there is no correlation between variables. The Q3 statistic is often used to test the residual correlation between items. If the absolute value of the Q3 statistic is less than 0.2, it indicates that the local independence assumption is satisfied. According to the calculation results, the absolute



residual correlations for all items in this study are less than 0.2, indicating that these items meet the local independence assumption.

### *Item fit test*

Item fit testing is used to evaluate how well test items align with the chosen Item Response Theory model, ensuring the quality and validity of the measurement. The S-$\chi^2$ statistic standardizes the traditional $\chi^2$ fit statistic, reducing the influence of sample size on the results and assessing item fit. A p-value greater than 0.05 indicates good fit. According to the calculations, the p-values for all items in this study are shown in Table 4, and all items exhibit good fit.

**Table 4**
*S-$\chi2$ and Parameters*

| item | p. S-$\chi2$ | discrimination | difficulty |
|---|---|---|---|
| 1 | 0.801 | 1.007 | 0.12 |
| 2 | 0.191 | 0.999 | 0.313 |
| 3 | 0.572 | 0.831 | -2.105 |
| 4 | 0.165 | 1.257 | 0.085 |
| 5 | 0.158 | 0.977 | -0.35 |
| 6 | 0.508 | 0.923 | -1.558 |
| 7 | 0.538 | 0.931 | 0.220 |
| 8 | 0.724 | 0.908 | 0.464 |
| 9 | 0.622 | 1.201 | -1.196 |
| 10 | 0.092 | 1.276 | -0.705 |
| 11 | 0.411 | 0.862 | 0.491 |
| 12 | 0.115 | 0.934 | 1.204 |
| 13 | 0.315 | 1.129 | 1.762 |
| 14 | 0.120 | 1.035 | 0.284 |
| 15 | 0.642 | 0.962 | -0.179 |
| 16 | 0.137 | 1.149 | 0.150 |
| 17 | 0.296 | 0.978 | -0.621 |
| 18 | 0.132 | 0.857 | 0.176 |
| 19 | 0.719 | 1.119 | -0.742 |
| 20 | 0.171 | 1.073 | 1.835 |
| 21 | 0.583 | 0.917 | -0.759 |



| 22 | 0.067 | 1.209 | 0.733 |
| 23 | 0.217 | 1.017 | 0.734 |
| 24 | 0.612 | 1.074 | -0.458 |
| 25 | 0.130 | 1.018 | -0.880 |
| 26 | 0.701 | 0.859 | -0.596 |
| 27 | 0.433 | 1.132 | -0.242 |
| 28 | 0.283 | 1.327 | 0.797 |
| 29 | 0.280 | 0.945 | 0.310 |
| 30 | 0.093 | 0.907 | -0.655 |
| 31 | 0.403 | 0.937 | -0.392 |
| 32 | 0.091 | 1.054 | -0.210 |

### *Test Parameters*

The calculation of difficulty and discrimination is crucial for evaluating test item quality. Discrimination measures how well an item distinguishes between different ability levels, with the slope of the item characteristic curve (ICC) representing its sensitivity. High-discrimination items have steeper slopes and effectively differentiate examinee abilities, while low-discrimination items do not. Discrimination values (parameter a) are categorized as follows: 0.01–0.34 (extremely low), 0.35–0.64 (low), 0.65–1.34 (moderate), 1.34–1.69 (high), and above 1.70 (very high). Items with a discrimination value below 0.64 are recommended for removal due to their limited contribution to ability measurement. The parameter calculations are shown in Table 4, and the discrimination of all items meets the requirements.

The total item characteristic curve, the item characteristic curves for each individual item and wright map are shown in the figures below.



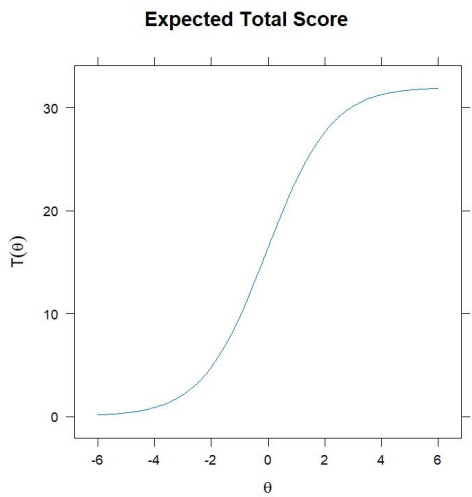

**Figure 2**

*Item Characteristic Curves*

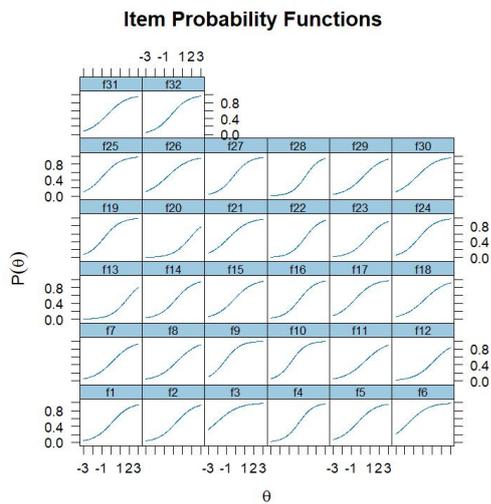

**Figure 3**

*Item Characteristic Curves for Each Item*

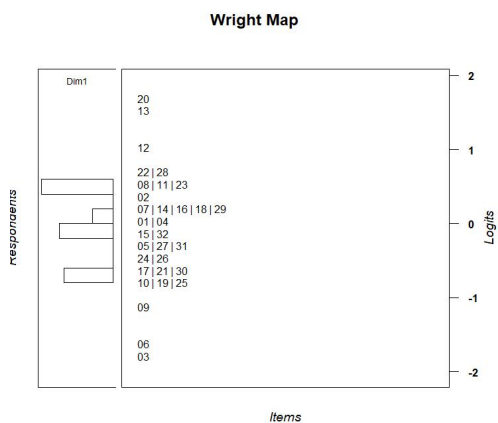



**Figure 4**
*Wright Map*

Based on these parameters, the difficulty of the items covers a wide range, which helps differentiate students' ability levels. The design of the item difficulties is reasonable, with relatively simple items (e.g., F3) as well as more challenging ones (e.g., F2). The distribution of the difficulty parameters shows that the test items are suitable for assessing students with different ability levels. The gradient from easy to difficult is well-balanced, making the test appropriate for evaluating a broad spectrum of mathematical proficiency or other related abilities.

***Variable Controls***

Given the uncertain variables involved in framing RQ2, the following controls are made.First, all recruited subjects have taken a complete one-semester course in linear algebra, and all had passed the final test at their schools (Hiebert, and Grouws 2007).Moreover, all the students studied linear algebra using the same textbook (Linear algebra, 2021), the very same textbook used earlier in the framework of RQ1. Second, all students are consistent in the timing of their PT sessions (Bray, 2014).Third, with regard to the quality of teaching, three professors of mathematics are used for teaching, who use the same teaching handouts and have been guided in their teaching by the researcher (Guill et al ., 2019; Zhang et al, 2021). And all three professors have PhDs in mathematics and have more than 10 years of experience teaching linear algebra.

***Description of formal tests***

The data set used by the research institute was obtained using convenience samples from a total of 60 students from different levels of schools in Zhejiang, Jiangsu, Hubei and



Shandong provinces. These students all found local tutoring institutions to study linear algebra to prepare for the mathematics test on the graduate school entrance exam.

These students were all third year undergraduate students, including 38 men and 22 women. Before tutoring, these students had basic knowledge of linear algebra and had completed a course on linear algebra. These students took a pretest before tutoring, engaged in three days of 24 lessons in total, and completed a posttest on the fourth day. The handouts used by the three teachers were uniformly arranged to ensure that the pretest and posttest topics would not be the same as those in the tutorial. The questionnaire was completed in the posttest and did not consume any of the testing time.

In the end, everyone was able to complete the pretest and posttest on time, and 60 questionnaires were collected. All statistical analyses were performed using SPSS 23. The IRT model and other data analyses that require drawing based on Python 3.3 were used.

The selection of questions for the pretest and posttest was based on a professional linear algebra problem set, accompanied by reference answers, and error-prone questions on knowledge points such as determinants, matrices, systems, and vectors that were provided in the previous text were used. In terms of difficulty measurement, an online questionnaire was used. The questions were measured online from 2021 to 2022 in the form of multiple-choice questionnaires, with 659 and 674 scores received on the pretest and posttest, respectively. The vast majority of the questionnaires came from educational institutions in Shandong, Jiangsu, Zhejiang, and other regions that used live class sessions as preclass training for students, and the difficulty was measured based on data from these initial sessions.

The difficulty curve of each question is shown in Figure 5. The corresponding pretest



and posttest IRT difficulty comparison chart is shown in Figure 6. The overall difficulty of

Test Two was higher than that of Test One, which is a proposition principle unanimously

obtained by the research team after discussion. Informal interviews with three mathematics

professors confirmed the belief that short-term training can achieve good results in solving

fixed problems for students, which is inconsistent with the original intention. To obtain a

more realistic understanding of students' abilities, in the posttest, the difficulty of each

corresponding question was increased, and therefore fixed modes could not be directly used

to solve it, which better aligned with the research questions of this study.

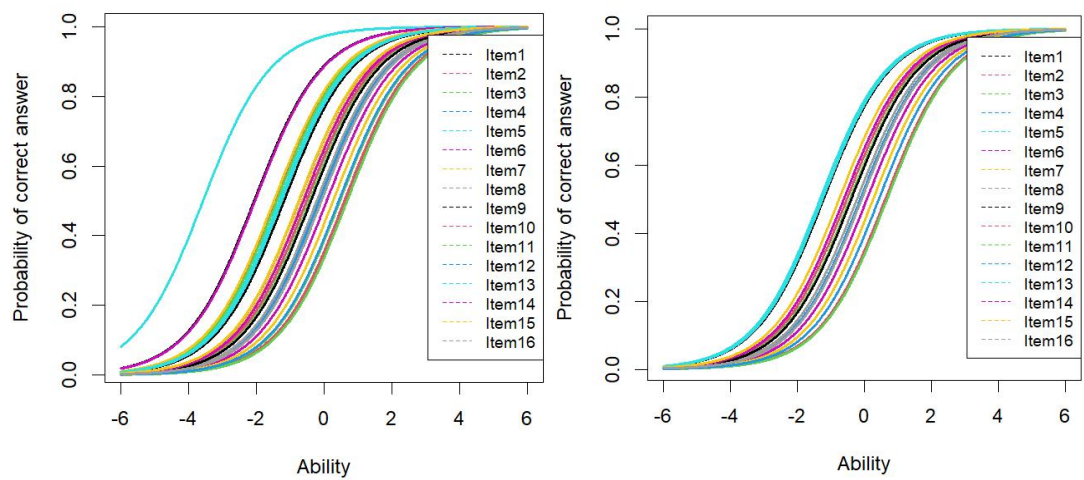

**Figure 5**
*Probability of Correct Answers for Pretest (Left) and Posttest (Right)*



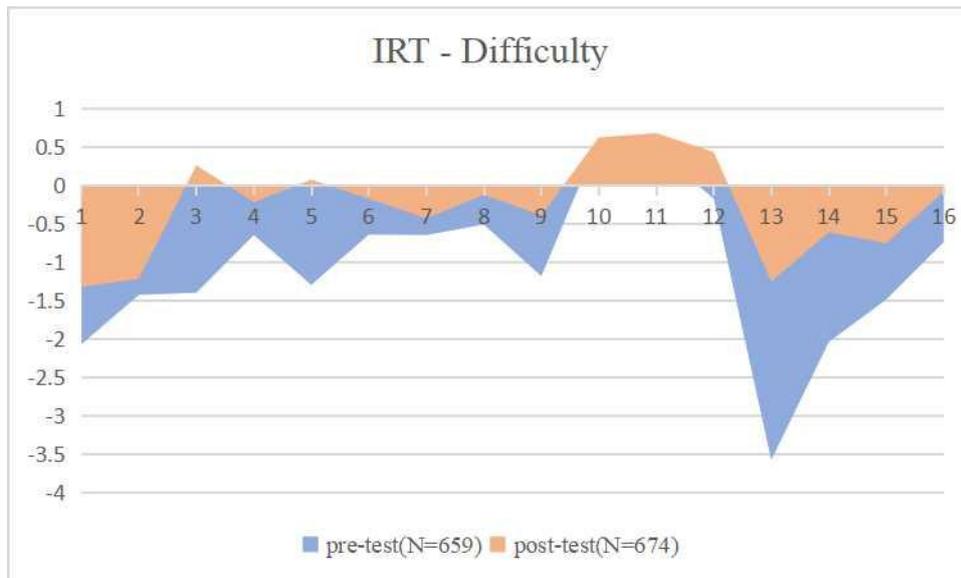

**Figure 6**
*IRT Difficulty Levels for Pretest and Posttest*

      Notably, in formal testing, offline paper and pencil tests are used, and the questions

are somewhat different from online questionnaire tests, which follows testing conventions in

the discipline of mathematics. A professor who participated in this study said, "We cannot

tolerate a single test being all multiple-choice questions. So the fill-in-the-blank section of the

online questionnaire was converted into a solution question, but the rating was still 0 or 1,

using the IRT model."

      This has caused some confusion in the grading of the question. For the

problem-solving section, if students have demonstrated the correct concept but there are

certain computational difficulties, it was determined after unanimous discussion to give them

a score of 1 because they showed mastery of this problem at the thinking level. For example,

it was recognized that the students in Figure 7 could achieve full marks in solving linear

equation systems because they mastered the methods of solving the system and the structure

of the solution was well represented. Although students' writing is not standardized and



meaningless constants such as $k_3$ may even appear, the basic solving process and ideas are already understood by students.

$$A = \begin{pmatrix} 1 & 1 & 0 & 0 & 1 \\ 1 & 1 & -1 & 0 & 0 \\ 0 & 0 & 1 & 1 & 1 \end{pmatrix} \rightarrow \begin{pmatrix} 1 & 1 & 0 & 0 & 1 \\ 0 & 0 & -1 & 0 & -1 \\ 0 & 0 & 0 & 1 & 0 \end{pmatrix}$$

$$\begin{cases} x_1 + x_2 + x_5 = 0 \\ x_2 = k_1 \\ -x_3 - x_5 = 0 \\ x_4 = 0 \\ x_5 = k_2 \end{cases} \Longrightarrow \begin{cases} x_1 = -k_1 - k_2 \\ x_2 = k_1 \\ x_3 = -k_2 \\ x_4 = 0 \\ x_5 = k_2 \end{cases}$$

$$(x_1, x_2, x_3, x_4, x_5)^T = (-1, 1, 0, 0, 0) k_1 + (-1, 0, -1, 0, 1) k_2$$

**Figure 7**
*An Answer to the Linear Systems Question*

**Multiple Linear Regression**

The dependent variable is the posttest latent traits in Table 5, and the remaining variables are independent variables. The formula is shown in (1).

$$y = \beta_0 + \beta_1 x_1 + \beta_2 x_2 + \cdots \beta_p x_p + \varepsilon \quad (1)$$

**Table 5**
*Assignment and Definitions of the Variables*

| Variables | Description |
|---|---|
| **Product variables: (Post-Pre)** | |
| Posttest latent traits | Standardized scores derived from the original scaling scores |
| Pretest latent traits | Standardized scores derived from the original scaling scores |
| **Process variables** | |
| Instructional content | _ |
| **Presage variables: Teacher qualifications** | |
| Student–teacher ratio | The number of students divided by the number of teachers in a class |
| **Presage variables: Time point** | |
| Study time | 1=None, 2=Less than 30 min, 3=30 to 60 min, 4=1 to 2 h, 5=More than 2 h |

**Context variables: Self-expectation of education**

| | |
|---|---|
| Interest in learning | 1=Very Low, 2=Low, 3=Medium, 4=High, 5=Very High |
| Self-evaluation | 1=Very Low ,2=Low, 3=Medium, 4=High, 5=Very High |

**Context variables: Family SES**

| | |
|---|---|
| Family's economic condition | 1= very poor, 2= poor, 3= medium, 4= rich, 5= very rich |
| Study table | If the student has an individual study table at home: 1= yes, 2= no |
| Books | Family's collection of books: 1 = few, 2=some, 3= moderate, 4= several, 5 = plenty |
| Phone or Pad | If the family has an iPhone or iPad: 1= yes, 2= no |

**Individual information about students: Demographic characteristics**

| | |
|---|---|
| Gender | 0=Female, 1=Male |
| Single child | 0=Non-single child, 1=Single child |
| Learning difficulty | Students' perception of level of difficulty in learning mathematics: 1= very, 2= a little, 3=not know, 4= not very, 5= none |
| Students' estimation | Students' estimated score in mathematics if full score is 100 |

The variable Instructional content does not appear in the final model because some of the variables have been controlled for in the variable control section. Not only that, but when answering the questionnaire, all of the questionnaires indicated that they had a separate desk and their own mobile phone when answering questions in the category of Family SES. Not only that, 55 students chose the middle when responding to family economic status and 51 students chose the middle when responding to family book collection, which made this question not well investigated, so the variable of family SES was omitted in the final model. This phenomenon is well explained by the fact that according to Guo (2021) there is a good



description of the characteristics of Chinese students, who are very modest and cautious.

To better measure self-efficacy, six questions were used in the questionnaire design, with a Cronbach coefficient of $\alpha = 0.97$. These four questions were added as a parameter to the multiple linear regression model to study students' latent traits. The formula for its definition is (2):

$$x_i = \frac{1}{6}(x_{i1} + x_{i2} + x_{i3} + x_{i4} + x_{i5} + x_{i6})\ (2)$$

**Results**

Through the IRT model, the latent traits of 60 students in the pretest and posttest were calculated, as shown in Figure 8. After the Wilcoxon signed rank test of relevant samples, the results are shown in Table 6. It can be concluded that the latent traits of students were significantly improved. From Figure 8, it can also be seen that for each individual student, their latent traits were improved through two tests.

**Table 6**
*Sign Rank Test Results of the Pretest and Posttest*

| $post - pre$ | $N$ | Rank average | $Z$ | $p$ |
|---|---|---|---|---|
| negative rank | 4 | 18.000 | | |
| positive rank | 55 | 30.870 | -6.139 | $<0.001$ |



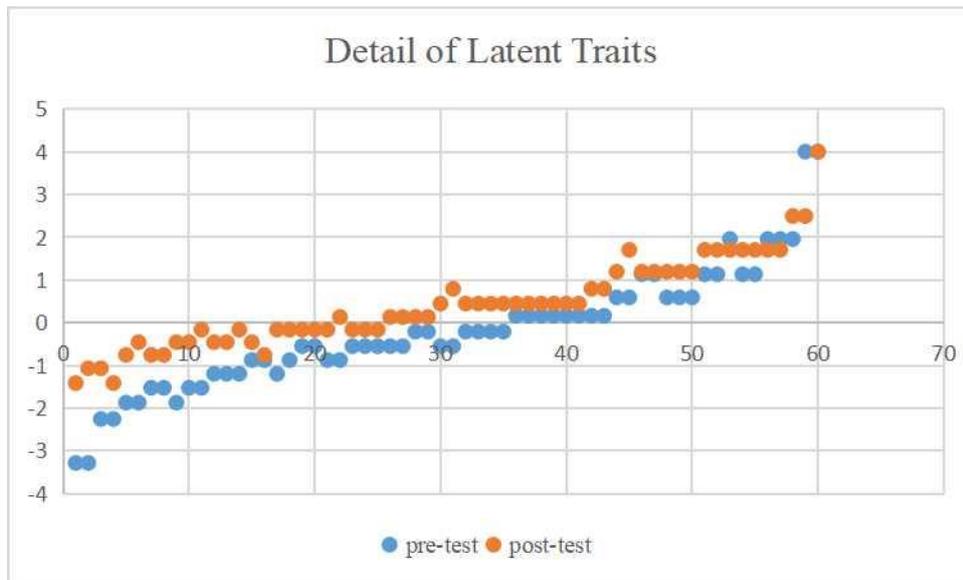

**Figure 8**
*Details of Latent Traits for 60 Students*

For the second research question, a multiple linear regression model was used to introduce various possible variables into the model using the forward method. The obtained regression parameters are shown in Table 7.

**Table 7**
*Normalized Linear Regression Model Coefficients*

| Category | $\beta_i$ | $t$ | Sgn. |
|---|---|---|---|
| Gender | -0.134 | -1.171 | 0.246 |
| Single child | 0.028 | 0.247 | 0.805 |
| Interest in learning | 0.081 | 0.701 | 0.486 |
| Student–teacher ratio | -0.205 | -1.809 | 0.076* |
| Study time | 2.009 | 3.058 | 0.005** |
| Self-evaluation | 0.791 | 2.993 | 0.003** |

$* \, p < 0.1 , ** \, p < 0.01$

Finally, a linear regression model was obtained using variable self-efficacy and variable learning duration ( $F = 10.163, p < 0.001$ ). The variable of teacher-student ratio is



not used here is because there is an insufficient amount of test data. As three teachers are teaching, there is no good quantitative support for the teacher-student ratio variable. Moreover, the latent traits of the pretest were also discarded; according to Table 6, it can be seen that a simple linear regression between the two may affect the position of other variables in the model, which is not in line with the original intention of research question 2. The same conditions apply to the actual situation: if the latent traits of the pretest are utilized, other variables will no longer have significance.

It is evident through Table 6 that the dependent variable (post-pre) is dependent on three significant independent variables namely Student-teacher ratio, Study time, and Self-evaluation. Due to the lack of instructors, we are unable to make statistical inferences about the variable of student-teacher ratio, but it can be assumed that, to some extent, the student-teacher ratio influences the improvement of student performance. And not just that, the smaller the student-teacher ratio, the more pronounced the improvement in student achievement. Not only that, but the time spent by the students on their studies is the right factor that influences the improvement of their scores, and it can be determined that the Studying for a short period of time can lead to a high level of knowledge acquisition, and this study is consistent with the results obtained by Zhang (2011). In addition, students can clearly perceive the improvement of their academic performance and this result is remarkable.

**Conclusion**

Based on the above data, the following conclusions can be drawn:

(1) The impact of PT on linear algebra in higher education



The study highlights the profound influence of PT on bolstering students' problem-solving acumen in the context of linear algebra in higher education. Interestingly, as students progress in their academic grades, the magnitude of the impact of appears to diminish or even retract. This suggests that coaching has a dual utility: for high achievers, it serves as a reinforcement mechanism to consolidate their strengths, and for students lagging behind, it offers a conduit to narrow the academic disparity with their peers (Bray & Kobakhidze, 2014; Zhang, 2011). While these findings diverge from established patterns observed in primary and secondary education (Ren et al., 2022), they resonate with Davies (2004), who posited the academic dividends of extended training in complex cognitive domains.

(2) Performance in Determinants:

Chinese students demonstrated remarkable adeptness in determinant-related problems. Across the two administered tests, their aptitude was evident, as they adeptly applied the result of $\det kA \neq k \det A$ with a notably low error rate. Such a clear understanding sharply contrasts with the misconceptions noted in Aygor and Ozdag's (2012) study. Furthermore, when juxtaposed with the responses from a substantial cohort of educators in Zimbabwe — specifically, the 116 teachers that were reported by Kazunga and Bansilal (2018) — the instructive methodologies of the three professors involved in this research stand out as particularly efficacious. A plausible explanation for such a stark distinction could be in the pedagogical approach adopted: all three professors unanimously emphasized the importance of determinants at the inception of linear algebra instruction. Their consistent reinforcement, especially through frequent testing on this specific topic, likely facilitated a robust conceptual



grasp among the students.

(3) Challenges with Matrix Multiplication:

Despite the influence of PT, students' skills in matrix multiplication remained relatively static, underscoring a gap in their foundational understanding. This observation aligns with Cook et al. (2018), who noted similar challenges when matrix multiplication was integrated with elementary matrices. The observed difficulty suggests that the mere introduction of matrix multiplication concepts can significantly challenge students, necessitating a re-evaluation of instructional methodologies.

(4) The linear equation conundrum:

The findings pertaining to linear equations corroborate established research paradigms. While students demonstrate an aptitude in discerning solutions to a system of equations, a conceptual gap emerges when attempting to ascertain whether a system's solution simultaneously satisfies individual equations within that system, a phenomenon previously observed by Panizza et al. (1999). Notably, students demonstrate a prevalent ability to compute accurate solutions to these systems of equations, even in the absence of PT interventions; their inherent capabilities allow them to swiftly find solutions. However, in grappling with the intricate nuances of the systems' solution structures, challenges arise. Following the pedagogical interventions of PT, although there have been gains in bridging this conceptual gap, they have been modest.

(5) Vector interpretation:

In the domain of vector interpretation, student performance was found to be suboptimal, aligning with the patterns elucidated by Ertekin et al. (2010). Such challenges



can be attributed to cognitive barriers that students encounter when navigating linear-associated problems within algebraic constructs. Interestingly, these impediments did not preclude the marked advancements witnessed through geometric interpretations, as shown in the work of Stewart et al. (2019). Through PT interventions, students not only deepened their grasp of geometric interpretations of linearity but also honed their capabilities to synergize linear equation problems with their linearly associated counterparts. This enhanced coordination, cultivated under the guidance of adept PT instructors, has significantly bolstered their proficiency in addressing analogous problems.

(6) Prism of self-efficacy:

In the context of self-efficacy, a discernible positive correlation emerges between students' self-awareness and their mathematical achievements as well as the duration of their learning. Intriguingly, this self-awareness appeared largely unaffected by demographic variables such as gender and single-child status. This observation resonates with Davies' (2004) assertion that sustained engagement in intricate cognitive training can yield academic dividends. Furthermore, this pattern parallels findings from a study on primary school attendees (Guo et al., 2021).

The post-PT session feedback was overwhelmingly affirmative. Students collectively expressed that this brief but intense tutoring phase catalyzed substantial advancements in their grasp of linear algebra. One student's feedback, which was representative of the group's sentiment, noted a breakthrough in understanding previously elusive concepts. Furthermore, the pedagogical experts echoed this optimism. The trio of professors involved in this intervention underscored the transformative potential of the PT sessions, positing that



judiciously designed PT is instrumental in refining students' conceptual clarity in linear algebra. This sentiment aligns with pedagogical insights from Canadian educational research (Davies, 2004) and corroborates findings presented by Bray and Kobakhidze (2014) and Zhang (2011).

**Learning Opinions and PT Suggestions**

Linear algebra presents a unique pedagogical challenge, characterized by its potential to mask superficial comprehension beneath the façade of algorithmic proficiency. While students may demonstrate aptitude in memorized proofs and known algorithms, such outward displays can often conceal deeper conceptual gaps. A student's ability to eloquently articulate the definition of a linear transformation, for instance, does not inherently guarantee their competence in identifying its application in varied contexts (Sandoval & Oktaç, 2021).

Furthermore, the caliber of instructor plays an instrumental role in shaping these student outcomes. A misinterpretation or oversimplification by the educator can cascade into profound misunderstandings for the student, distracting them from the foundational principles and even impeding their ability to grasp core content (Kazunga & Bansilal, 2018). While the faculty of larger, well-established institutions is often well-equipped to navigate the intricacies of linear algebra, feedback from various institutional leaders has underscored a concerning trend among smaller institutions with limited staff: a discernible deficit in robust, comprehensive instruction. Such gaps, if unaddressed, risk leaving students with a fragmented understanding of pivotal concepts.

In teaching mathematical concepts to students, it is crucial to avoid overwhelming them with an abundance of intricate details. Instead, educators should focus on key



overarching strategies that promote deeper conceptual understanding. The suggestions for promoting these strategies are as follows (Oktaç, 2018, 2019; Trigueros, 2018):

(1) Employ tasks that are unfamiliar to students, thereby encouraging a conceptual approach rather than a mere reliance on algorithms.

(2) Upon encountering mistakes in students' problem-solving attempts, it is beneficial to present the same problem from diverse perspectives. This facilitates a comparison between contrasting interpretations, enabling students to identify and rectify their errors.

(3) Craft mathematical scenarios tailored for the application of cognitive processes, such as internalization and encapsulation.

(4) Incorporate a variety of representations—including linguistic, matrix-based, graphical, algebraic, and numerical ones. Crucially, educators should guide students in discerning the connections between these representations, fostering a seamless transition among them.

(5) Recognize that students require adequate time to synthesize knowledge. It is during these reflective intervals that they may experience moments of epiphany or shifts in conceptual understanding, enabling them to tackle problems with renewed creativity and insight.

In informal interviews with three esteemed teaching professors, a shared sentiment emerged: an aversion to excessively procedural teaching approaches. Such approaches often lead to the mere imprinting of algorithms in students' minds, neglecting the intrinsic essence of the subject matter. This procedural emphasis impedes students' ability to seamlessly



integrate new concepts with foundational knowledge upon transitioning between topics. This observation aligns with the theoretical framework posited in the concept of Tall's 3 worlds (Beltrán-Meneu et al., 2016). While the professors did not articulate their insights in theoretical parlance, their pedagogical practices consistently mirror these findings.

For students who have skillfully internalized problem-solving concepts, the utility of procedural teaching (PT) becomes less pronounced. Moreover, high-achieving students tend to derive diminishing marginal benefits from PT. However, for students who grapple with problem solving and possess a robust sense of self-efficacy, PT remains indispensable. Its application not only bolsters their confidence in dealing with linear algebra but also markedly enhances their mathematical performance.

## Limitations and Future Work

The current study is characterized by several limitations:

(1) Measurement of test difficulty

The assessment of test difficulty was contingent upon online questionnaires, a modality that gained traction during the COVID-19 pandemic. The method's inherent constraints make it challenging to ascertain whether participants utilized external aids during the test. Nevertheless, the difficulty metrics for both the pretest and posttest showed remarkable consistency.

(2) Sample size considerations:

While the offline testing sample size exceeded that in many higher education studies, it falls short when compared to large-scale studies, casting doubts on data reliability. This limitation, which is pervasive in higher education research, stems from the granular



specialization of subjects and consequent restricted sample sizes. Additionally, our participants, representing the academic years 2022 and 2023, might exhibit varied latent traits. This variation could be attributed to disparities in college entrance examination scores in the Chinese mainland across institutions and years.

(3) Intricacies of PT:

PT often adopts a more informal, less structured, and voluntary methodology, diverging from traditional institutional norms. PT environments tend to be less conventional than mainstream educational settings, as evidenced by the prevalence of PT sessions held in personal spaces such as student or teacher residences (Zhang et al., 2021). Such nuances pose inherent challenges to the research, complicating data collection processes. Furthermore, while PT instruction does not starkly deviate from traditional pedagogical methods, it possesses unique attributes (Sandoval, & Oktaç,, 2021). These distinct characteristics may not be readily discernible using evaluative instruments designed for mainstream educational contexts.

Building upon the findings of the current research, subsequent studies should further explore both the depth of mathematical course comprehension and the self-perception of students in broader science and engineering domains. To do so, the research framework could be extended to disciplines such as calculus and probability theory. Such endeavors would be pivotal in ascertaining whether PT offers tangible pedagogical benefits in tertiary education. A critical inquiry would be to discern whether PT, rather than overburdening learners, as has been observed at the primary and secondary levels, truly promotes mathematical proficiency in a higher education context.



# References

Arnon, I., Cottrill, J., Dubinsky, E., Oktaç, A., Fuentes, S. R., Trigueros, M., & Weller, K. (2014). *APOS theory—a framework for research and curriculum development in mathematics education*. Springer.

Axler, S. (1997). *Linear algebra done right*. Springer Science & Business Media.

Aygor, N., & Ozdag, H. (2012). Misconceptions in linear algebra: The case of undergraduate students. *Procedia - Social and Behavioral Sciences, 46*, 2989-2994. https://doi.org/10.1016/j.sbspro.2012.05.602

Bagley, S., & Rabin, J. M. (2015). Students' use of computational thinking in linear algebra. *International Journal of Research in Undergraduate Mathematics Education, 2*(1), 83-104. https://doi.org/10.1007/s40753-015-0022-x

Baker, D. P., Motoko, A., LeTendre, G. K., & Wiseman, A. W. (2001). Worldwide shadow education: Outside-school learning, institutional quality of schooling, and cross-national mathematics achievement. *Educational Evaluation and Policy Analysis, 23*(1), 1-17. https://doi.org/10.3102/01623737023001001

Beltrán-Meneu, M. J., Murillo-Arcila, M., & Albarracín, L. (2016). Emphasizing visualization and physical applications in the study of eigenvectors and eigenvalues. *Teaching Mathematics and its Applications, 36*(3), hrw018. https://doi.org/10.1093/teamat/hrw018

Bray, M. (2014). The impact of shadow education on student academic achievement: Why the research is inconclusive and what can be done about it. *Asia Pacific Education Review, 15*(3), 381-389. https://doi.org/10.1007/s12564-014-9326-9

Bray, M., & Kobakhidze, M. N. (2014). Measurement issues in research on shadow education: Challenges and pitfalls encountered in TIMSS and PISA. *Comparative Education Review, 58*(4), 590-620. https://doi.org/10.1086/677907

Bray, M., & Kwok, P. (2003). Demand for private supplementary tutoring: Conceptual considerations, and socio-economic patterns in Hong Kong. *Economics of Education Review, 22*(6), 611-620. https://doi.org/10.1016/s0272-7757(03)00032-3

Carroll, J. B. (1970). *A model of school learning.* (Vol. 64, pp. 723–733). Cambridge, Mass.

(Original work published 1963)

Cook, J. P., Zazkis, D., & Estrup, A. (2018). Rationale for matrix multiplication in linear algebra textbooks. In S. Stewart, C. Andrews-Larson, A. Berman, & M. Zandieh (Eds.), *Challenges and strategies in teaching linear algebra ICME-13 monographs* (pp. 103-125). Springer International Publishing.

Davies, S. (2004). School choice by default? Understanding the demand for private tutoring in Canada. *American Journal of Education, 110*(3), 233-255. https://doi.org/10.1086/383073

DeMars, C. (2010). *Item response theory*. Oxford University Press.

Department of Mathematics, Shanghai Jiao Tong University. (2021). Linear algebra (3rd ed.). Science Publishing House. (Original work published 2000) ,ISBN



978-7-03-041199-0.

DeVries, D., & Arnon, I. (2004). Solution—What does it mean? Helping linear algebra students develop the concept while improving research tools. In M. J. Høines & A. B. Fuglestad (Eds.), *Proceedings of the 28th conference of the international group for the psychology of mathematics education* (Vol. 2, pp. 55–62). Bergen University College.

Dias, M. A., & Artigue, M. (1995). Articulation problems between different systems of symbolic representations in linear algebra. In L. Meira & D. Carraher (Eds.), *Proceedings of the 19th PME Conference* (Vol. 2, pp. 34–41). Atual Editora.

Dogan-Dunlap, H. (2010). Linear algebra students' modes of reasoning: Geometric representations. *Linear Algebra and its Applications, 432*(8), 2141-2159. https://doi.org/10.1016/j.laa.2009.08.037

Dogan-Dunlap, H. (2018). Mental schemes of: Linear algebra visual constructs. In S. Stewart, C. Andrews-Larson, A. Berman, & M. Zandieh (Eds.), *Challenges and strategies teaching linear algebra* (pp. 219-239). Springer.

Donevska-Todorova, A. (2014). Three modes of description and thinking of linear algebra concepts at upper secondary education. *Beiträge zum Mathematikunterricht, 48*(1), 305–308.

Dorier, J.-L. (1995). Meta level in the teaching of unifying and generalizing concepts in mathematics. *Educational Studies in Mathematics, 29*(2), 175-197. https://doi.org/10.1007/bf01274212

Dorier, J.-L., & Sierpinska, A. (2001). Research into the teaching and learning of linear algebra. In D. Holton & M. Artigue (Eds.), *The teaching and learning of mathematics at university level: An ICMI study* (pp. 255-273). Springer.

Dunkin, M. J., & Biddle, B. J. (1974). The study of teaching. Holt, Rinehart & Winston.

Embretson, S. E., & Reise, S. P. (2013). *Item response theory*. Psychology Press.

Ertekin, E., Solak, S., & Yazici, E. (2010). The effects of formalism on teacher trainees' algebraic and geometric interpretation of the notions of linear dependency/independency. *International Journal of Mathematical Education in Science and Technology, 41*(8), 1015-1035. https://doi.org/10.1080/0020739x.2010.500689

Ferryansyah, Widyawati, E., & Rahayu, S. W. (2018). The analysis of students' difficulty in learning linear algebra. *Journal of Physics: Conference Series, 1028*, 012152. https://doi.org/10.1088/1742-6596/1028/1/012152

Guill, K., Lüdtke, O. and Köller, O. (2020), Assessing the instructional quality of private tutoring and its effects on student outcomes: Analyses from the German National Educational Panel Study. Br J Educ Psychol, 90: 282-300. https://doi.org/10.1111/bjep.12281

Guo, M., Hu, X., & Leung, F. K. S. (2021). Culture, goal orientations, and mathematics achievement among Chinese students. *International Journal of Science and Mathematics Education, 20*(6), 1225-1245. https://doi.org/10.1007/s10763-021-10202-0

Hambleton, R. K., Swaminathan, H., & Rogers, H. J. (1991). *Fundamentals of item response*



*theory (Vol. 2)*. Sage.

Harel, G. (1987). Variations in linear algebra content presentations. *For the Learning of Mathematics, 7*(3), 29–32.

He, Y., Zhang, Y., Ma, X., & Wang, L. (2021). Does private supplementary tutoring matter? The effect of private supplementary tutoring on mathematics achievement. *International Journal of Educational Development, 84*, 102402. https://doi.org/10.1016/j.ijedudev.2021.102402

Hiebert, J., & Grouws, D. A. (2007). The effects of classroom mathematics teaching on students' learning. Second handbook of research on mathematics teaching and learning, 1(1), 371-404.

Hillel, J. (2000). Modes of description and the problem of representation in linear algebra. In J.-L. Dorier (Ed.), *On the teaching of linear algebra* (pp. 191-207). Springer.

Kazunga, C., & Bansilal, S. (2018). Misconceptions about determinants. In S. Stewart, C. Andrews-Larson, A. Berman, & M. Zandieh (Eds.), *Challenges and strategies in teaching linear algebra ICME-13 monographs* (pp. 127-145). Springer.

Lang, S. (2012). *Introduction to linear algebra*. Springer Science & Business Media.

Liu, J., & Bray, M. (2018). Evolving micro-level processes of demand for private supplementary tutoring: Patterns and implications at primary and lower secondary levels in China. *Educational Studies, 46*(2), 170-187. https://doi.org/10.1080/03055698.2018.1555452

Lord, F. M. (2012). *Applications of item response theory to practical testing problems*. Routledge.

Mischo, C., & Haag, L. (2002). Expansion and effectiveness of private tutoring. *European Journal of Psychology of Education, 17*(3), 263-273. https://doi.org/10.1007/bf03173536

Oktaç, A. (2018). Conceptions about system of linear equations and solution. In S. Stewart, C. Andrews-Larson, A. Berman, & M. Zandieh (Eds.), *Challenges and Strategies in Teaching Linear Algebra ICME-13 Monographs* (pp. 71-101). Springer.

Oktaç, A. (2019). Mental constructions in linear algebra. *ZDM, 51*(7), 1043-1054. https://doi.org/10.1007/s11858-019-01037-9

Panizza, M., Sadovsky, P., & Sessa, C. (1999). La ecuación lineal con dos variables : Entre la unicidad y el infinito. *Enseñanza de las Ciencias Revista de investigación y experiencias didácticas, 17*(3), 453-461. https://doi.org/10.5565/rev/ensciencias.4073

Ren, P., Dou, Z., Wang, X., Li, S., & Wang, L. (2022). Think twice before seeking private supplementary tutoring in mathematics: A data set from China questioned its effectiveness. *The Asia-Pacific Education Researcher, 32*(3), 429-437. https://doi.org/10.1007/s40299-022-00664-3

Samejima, F. (1997). Graded response model. In W. J. van der Linden & R. K. Hambleton (Eds.), *Handbook of modern item response theory* (pp. 85-100). Springer.

Sandoval, O. R., & Oktaç, A. (2021). Modelos intuitivos sobreel concepto de transformación lineal. In *Actes du colloque hom-mage à Michèle Artigue*. Université Paris Diderot.

Sierpinska, A. (2000). On some aspects of students' thinking in linear algebra. In J.-L. Dorier (Ed.), *On the teaching of linear algebra* (pp. 209-246). Kluwer.




Sierpinska, A., Dreyfus, T., & Hillel, J. (1999). Evaluation of a teaching design in linear algebra: The case of linear transformations. *Recherches en Didactique des Mathématiques, 19*(1), 7–40.

Stewart, S., Andrews-Larson, C., & Zandieh, M. (2019). Linear algebra teaching and learning: Themes from recent research and evolving research priorities. *ZDM, 51*(7), 1017-1030. https://doi.org/10.1007/s11858-019-01104-1

Strang, G. (2022). *Introduction to linear algebra*. Wellesley-Cambridge Press.

Trigueros, M. (2018). Learning linear algebra using models and conceptual activities. In S. Stewart, C. Andrews-Larson, A. Berman, & M. Zandieh (Eds.), *Challenges and strategies in teaching linear algebra ICME-13 monographs* (pp. 29-50). Springer International Publishing.

Wainer, H., Dorans, N. J., Flaugher, R., Green, B. F., & Mislevy, R. J. (2000). *Computerized adaptive testing: A primer*. Routledge.

Wang, L., & Guo, K. (2017). Shadow education of mathematics in China. In Y. Cao & F. K. S. Leung (Eds.), *The 21st century mathematics education in China* (pp. 93-103). Springer.

Wang, L., Gong, X., Pei, C.(2019). 数学教育视角下的影子教育研究的意义与未来方向 [Shadow education research with mathematics education: its value and future research area]. J. Math. Educ. 28 (1), 79–82.

Wawro, M., Zandieh, M., Rasmussen, C., & Andrews-Larson, C. (2013). *Inquiry oriented linear algebra: Course materials*. Creative Commons Attribution-NonCommercial-ShareAlike 4.0 International License.

Zeng, X., & Zhou, H. (2012). 北京市四、八年级学生课后补习的代价与收益 [A descriptive analysis on after-school tutoring in Beijing: Its costs and benfits]. *Journal of Educational Studies, 8*(6), 103-109.

Zhang, Y. (2011). *The determinants of national college entrance exam performance in [People's Republic of] China—with an analysis of private tutoring* [Unpublished Doctoral Dissertation]. Columbia University.

Zhang, Y., Dang, Y., He, Y., Ma, X., & Wang, L. (2021). Is private supplementary tutoring effective? A longitudinally detailed analysis of private tutoring quality in China. *Asia Pacific Education Review, 22*(2), 239-259. https://doi.org/10.1007/s12564-021-09671-3

Zumbo, B. D. (2007). Three generations of DIF analyses: Considering where it has been, where it is now, and where it is going. *Language Assessment Quarterly, 4*(2), 223-233. https://doi.org/10.1080/15434300701375832